\documentstyle[11pt]{article}
\begin{document}
\begin{center}\Large\bf{RANDOM SPARSE UNARY PREDICATES}\end{center}

\section{Introduction.}  Let $n$ be a positive integer, $0\leq p\leq 1$.
The random unary predicate $U_{n,p}$ is a probability space over
predicates $U$ on $[n]=\{1,\ldots,n\}$ with the probabilities determined
by
\[ \Pr[U(x)]=p,\mbox{    }1\leq x\leq n  \]
and the events $U(x)$ being mutually independent over $1\leq x\leq n$.  
Informally, we think of flipping a coin for each $x$ to determine if
$U(x)$ holds, the coin coming up ``heads'' with probability $p$.
We shall examine the first order language 
$<[n],\leq,U>$ with equality, a unary 
predicate $U$ and a binary predicate $\leq$.  Examples of sentences
in this language are:
\[ A:\exists_xU(x)  \]
\[ B:\exists_xU(x)\wedge\forall_y\neg y<x  \]
\[ C:\exists_{x,y}U(x)\wedge U(y)\wedge \forall_z\neg[x<z\wedge z<y]  \]
($>,\geq,<$ are natuarally definable from $\leq$ and equality.)
For any such sentence $S$ we have the probability
\[ \Pr[U_{n,p}\models S]  \]
While the use of unary predicates is natural for logicians there are
two other equivalent formulations  that will prove useful.
We may think of $U$ as a subset of $[n]$ and speak about
$i\in U$ rather than $U(i)$.  Second we may associate with
$U$ a sequence of zeroes and ones where the $i$-th term
is one if $U(i)$ and zero if $\neg U(i)$.  Thus we may
talk of starting at $i$ and going to   the next one.  We
shall use all three formulations interchangably.
\par Ehrenfeucht [??] showed that for any constant $p$
and any sentence $S$ in this language
\[ \lim_{n\rightarrow\infty} \Pr[U_{n,p}\models S]  \]
exists. In the case of sentences $A$ and $C$ the limiting probability
is one.  But  sentence $B$ effectively states  $1\in U$, hence
its limiting probability is $p$.  We get around these edge effects
with a new language, consisting of equality, a unary predicate $U$,
and a ternary predicate $C$.  We consider $C$ as a built in
predicate on $[n]$ with $C(x,y,z)$ holding if and only if either
$x<y<z$ or $y<z<x$ or $z<x<y$.  Thinking of $[n]$ as a cycle, with
$1$ coming directly after $n$, $C(x,y,z)$ is the relation that
$x$ to $y$ to $z$ goes in a clockwise direction.  For any 
sentence $S$ in this new language we can again define 
$\Pr[U_{n,p}\models S]$ only in this case Ehrenfeucht's results 
give a Zero-One Law: for any constant $p$ and sentence $S$
\[ \lim_{n\rightarrow\infty}\Pr[U_{n,p}\models S] = 0\mbox{ or }1 \]
We shall call the first language the {\em linear} language and the
second language the {\em circular} language. As a general guide, the
circular language will tend to Zero-One Laws while the linear
language, because of edge effects, will tend to limit laws.
\par We shall not restrict ourselves to $p$ constant but rather
consider $p=p(n)$ as a function of $n$.  We have in mind the
``Evolution of Random Graphs'' as first developed by Erd\H{o}s
and R\'enyi.  Here as $p=p(n)$ evolves from zero to one the 
unary predicate evolves from holding for no $x$ to holding for
all $x$.  Analogously (but without formal definition) we have
{\em threshold functions} for various properties.  For example,
$p(n)=n^{-1}$ is a threshold property for $A$.  When $p(n)\ll  
n^{-1}$ almost surely $A$ fails while when $p(n)\gg n^{-1}$
almost surely $A$ holds.  In Shelah,Spencer [??] we showed that
when $p=n^{-\alpha}$ with $\alpha\in(0,1)$, irrational then a 
Zero-One Law held for the random graph $G(n,p)$ and in \L uczak,
Spencer [??] we found a near characterization of those $p=p(n)$
for which the Zero-One Law held.  The situation with random unary
predicates turns out to be somewhat simpler. Let us say $p=p(n)$
satisfies the Zero-One Law for  circular unary predicates if for
every sentence $S$ in the circular language
\[ \lim_{n\rightarrow\infty}\Pr[U_{n,p(n)}\models S]=0\mbox{ or }1 \]
Here is our main result.

{\bf Theorem 1.} Let $p=p(n)$ be such that $p(n)\in[0,1]$ for all
$n$ and either
\[ p(n)\ll n^{-1}  \]
or for some positive integer $k$
\[ n^{-\frac{1}{k}}\ll p(n)\ll n^{-\frac{1}{k+1}}  \]
or for all $\epsilon > 0$
\[ n^{-\epsilon}\ll p(n) \mbox{ and } n^{-\epsilon}\ll 1-p(n) \]
or for some positive integer $k$
\[ n^{-\frac{1}{k}}\ll 1-p(n)\ll n^{-\frac{1}{k+1}}  \]
or
\[ 1-p(n)\ll n^{-1}  \]
Then $p(n)$ satisfies the Zero-One Law for circular unary predicates.
Inversely if $p(n)$ falls into none of the above categories then
it does not satisfy the Zero-One Law for circular unary predicates.

\par The inverse part is relatively simple.  Let $A_k$ be the sentence
that there exist $k$ consecutive elements $x_1,\ldots,x_k\in U$. 
($x,y$ are consecutive if for no $z$ is $C(x,z,y)$. For $k=2$ this is example C.
)  Then $\Pr[A_k]$ is (for a given $n$) a monotone function
of $p$.  When $p(n)\sim cn^{-1/k}$ and $c$ a positive constant the
probability $\Pr[A_k]$ approaches a limit strictly between zero and
one.  (Roughly speaking, $n^{-1/k}$ is a threshold function for $A_k$.)
Thus for $p(n)$ to satisfy the Zero-One law we must have $p(n)\ll n^
{-1/k}$ or $p(n)\gg n^{-1/k}$.  Further (replacing $U$ with $\neg U$),
the same holds with $p(n)$ replaced by $1-p(n)$.  For $p(n)$ to fall
between these cracks it must be in one of the above five categories.

\vspace{1cm}

\par {\em Remark.}  Dolan [??] has shown that $p(n)$
satisfies the Zero-One Law for
{\em linear} unary predicates if and only if $p(n)\ll n^{-1}$ or
$n^{-1}\ll p(n)\ll n^{-1/2}$ or $1-p(n)\ll n^{-1}$ or $n^{-1}\ll
1-p(n)\ll n^{-1/2}$.  For $n^{-1/2}\ll p(n)=o(1)$ he considered
the following property:
\[ D: \exists_xU(x)\wedge [U(x+1)\vee U(x+2)]\wedge \neg\exists_y[
U(y)\wedge [U(y+1)\vee U(y+2)]\wedge y<x] \wedge U(x+1)  \]
(Addition is {\em not} in our language but we write $x+1$ as
shorthand for that $z$ for which $x<z$ but there is no $w$
with $x<w<z$.)  In  our zero-one formulation
$D$ basically states that the first
time we have $11$ comes before the first time
we have $101$.  This actually has limiting probability
$.5$.  This example illustrates that limiting probability for linear
unary predicates can depend on edge effects and not just edge effects
looking at $U$ on a fixed size set $1,\ldots,k$ or $n,n-1,\ldots,n-k$.
We defer our results for linear unary predicates to section 4. 

\vspace{1cm}

When $p(n)\ll n^{-1}$ the Zero-One Law is trivially satisfies since
almost surely there is no $x$ for which $U(x)$.  Also, if $p(n)$
satisfies the Zero-One Law so does $1-p(n)$. Suppose $p=p(n)$
satisfies $p(n)\gg n^{-\epsilon}$ and $1-p(n)\gg n^{-\epsilon}$ for
all $\epsilon>0$.
We show in a section 3
that for every $t$ there is a sequence $A_1\cdots A_R$ with the property
that for any sentence $A$ of quantifier depth $t$ either all models
$<[u],C,U>$ that contain $A_1\cdots A_R$ as a subsequence satisfy $A$
or no such models satisfy $A$.  ($<[u],C,U>$ contains $A_1\cdots A_R$
as a subsequence
if for some $1\leq j\leq u$ for all $1\leq i\leq R$ we have $U(j+i)$
if and only if $x_i=1$, with $j+i$ defined modulo $u$.) For $p(n)$
in this range $<[u],C,U>$ almost surely contains any such fixed
sequence $A_1\cdots A_R$ as a subsequence and hence the Zero-One
Law is satisfied.
This leaves us with only one case in Theorem 1, and that will
be the object of the next section.
\section{The Main Case.}
Here we let $k$ be a positive integer and assume
\[ n^{-\frac{1}{k}} \ll p(n) \ll n^{-\frac{1}{k+1}}  \]
Our object is to show that $p(n)$ satisfies the Zero-One Law for
circular unary predicates. We shall let $t$ be a fixed, though
arbitrary large, positive integer.  We shall examine the equivalence
class under the $t$-move Ehrenfeucht game of the circular model.
For the most part, however, we shall examine linear models.
\par We define (as Ehrenfeucht did) an equivalence class on
models $M=<n,\leq,U>$, two models $M,M'$ being equivalent if
they satisfy the same depth $t$ sentences or, equivalently,
if the $t$-move Ehrenfeucht game on $M,M'$ is won by the
``Duplicator''.  The addition of models (with $M$ on $[n]$,
$M'$ on $[n']$ we define $M+M'$ on $[n+n']$) yields an
addition of equivalence classes.  We shall denote the 
equivalence classes by $x,y,\ldots$ and the sum by
$x+y$.  Results from the beautiful theory of these classes
are given in Section 3.
\par Let us consider a random unary function $U$ defined
on {\em all} positive integers
$1,2,\ldots$ and with $\Pr[U(i)]=p$ for all $i$, these
events mutually independent. (In the end only the
values of $U(i)$ for $1\leq i\leq n$ will ``count'' but
allowing $U$ to be defined over all positive integers
allows for a ``fictitious play'' that shall simplify
the analysis.) Now for any starting point
$i$ examine $i,i+1,\ldots$ until reaching the first $j$
(perhaps $i$ itself) for which $U(j)$.  Call $[i,j]$ the
1-interval of $i$. (With probability one there will
be such a $j$; fictitious play allows us to postpone
the analysis of those negligible cases when no $j$ is
found before $j>n$.) What are the possible Ehrenfeucht
values of $<[i,j],\leq,U>$?  The model must have a series
of zeroes (i.e., $\neg U$) followed by one one (i.e., $U$).
There is an $s$ ($s=3^t$ will do) so that all such models
with at least $s$ zeroes have the same Ehrenfeuct value.
We can write these values as $a_1,\ldots a_s$ and $b$
($a_i$ having $i-1$ zeroes, $b$ having $s$ zeroes).  Call
this value the {\em 1-value} of $i$.  The probability of the 
1-value being any particular $a_i$ is $\sim p$ while the
probability of it being $b$ is $\sim 1$.  (All asymptotics
are as $p\rightarrow 0$.)  We let $E_1$ denote this set of
possible 1-values and we split $E_1=P_1\cup T_1$
with $P_1=\{b\}$ and $T_1=\{a_1,\ldots,a_s\}$.  The 1-values
in $T_1$ we call 1-transient, the 1-value in $P_1$ we call
1-persistent.
\par Now (with an eye toward induction) we define the
2-interval of $i=i_0$.  Take the 1-interval of $i$, say
$[i_0,i_1)$.  Then take the 1-interval of $i_1$, say
$[i_1,i_2)$.  Continue until reaching a 1-interval
$[i_u,i_{u+1})$ whose 1-value is 1-transient.  (Of
course, this could happen with the very first interval.)
We call $[i,i_{u+1})$ the 2-interval of $i$.  Now we 
describe the possible 2-values for this 2-interval.
In terms of Ehrenfeucht value we
can write the interval as $b+b+\ldots+b+a_i$ where there
are $u$ (possibly zero) $b$'s.  Any $b+\ldots+b$ with
at least $s$ addends $b$ has (see \S 3.4) the same value, call it
$B$.  Let $jb$ denote the sum of $j$ $b$'s. We define
the transient 2-values $T_2$ as those of the form 
$jb+a_i$ with $0\leq j<s$ and the persistent 2-values  $P_2$
as those of the form $B+a_i$.  For example, let $t=5$ and
$s=3^5=243$.  Then $i$ has 2-value $6b+a_{22}$ if, starting
at $i$, six times there are at least 243 zeroes before a 
one and after the sixth one there are 21 zeroes and then a
one.  The 2-value is $B+a_5$ if at least 243 times there are
at least 243 zeroes before the next one and the first time
two ones appear less than 243 apart they are exactly 5 
apart.  What are the probabilities for $i=i_0$ having any
particular 2-value?  The first 1-interval $[i_0,i_1)$ has
distribution for 1-value as previously discussed: $\sim p$
for each $a_i$ and $\sim 1$ for $b$.  Having determined
the first 1-interval the values starting at $i_1$ have not
yet been examined.  Hence the 1-value of the second 1-interval
will be independent of the 1-value of the first and, in general,
the sequence of 1-values will be of mutually independent values.
Then the transient 2-values $jb+a_i$ each have probability $\sim p$
while the persistent 2-values $B+a_i$ will each have probability
$\frac{1}{s}+o(1)$. We let $P_2$ denote the set of persistent
2-values, $T_2$ the set of transient 2-values and $E_2=P_2\cup 
T_2$ the set of 2-values.
\par The 3-value will contain all the notions of the general case.
Beginning at $i=i_0$ take its 2-interval $[i_0,i_1)$.  Then
take successive 2-intervals $[i_1,i_2),\ldots,[i_{u-1},i_u)$ 
until reaching an interval $[i_u,i_{u+1})$
whose 2-value is transient.  The 3-interval for $i$
is then $[i,i_{u+1})$.  Let $x_1,\ldots,x_u,y_{u+1}$ be the
2-values for the successive intervals.
\relax Fromthe procedure all $x_i\in P_2$  while $y_{u+1}\in T_2$.
Now consider (see \S 3.1)the Ehrenfeucht equivalence classes
(again with respect to a $t$-move game)  over $\Sigma P_2$.
($\Sigma A$ is the set of strings over alphabet $A$.) Let $\alpha$
be the equivalence class for the string $x_1\cdots x_u$, then
the 3-value of $i$ is defined as
the pair $\beta=\alpha y_{u+1}$.  We let $E_3$ be
the set of all such pairs and we call $\beta$ persistent
(and place it in $P_3$) if $\alpha$ is a persistent state 
(as defined in \S 3.2)in $\Sigma P_2$;
otherwise we call $\beta$ transient and
place it in $T_3$. If $x_1\cdots x_u$ and $x_1'\cdots x_{u'}'$ are
equivalent as strings in $P_2$ then $x_1+\ldots+x_u$ and
$x_1'+\ldots+x_{u'}'$ have (as shown in \S 3.4)
the same Ehrenfeucht value.  So the
3-value of $i$ determines the Ehrenfeucht value of the 3-interval
of $i$ though possibly it has more information.
What are the  probabilities for the 3-value
of $i$?  Again we get a string of 2-values $z_1z_2\ldots$ whose
values are mutually independent and we stop when we hit a transient
2-value.   We shall see (in the course of the full induction argument)
that the probability of having 3-value $\beta$ is $\sim c_{\beta}$
for persistent $\beta$ and $\sim c_{\beta}p$ for transient $\beta$.
\par Now let us define $k$-interval and $k$-value, including the
split into persistent and transient $k$-values
by induction on $k$.  Suppose $E_k,P_k,T_k$ have been defined.  Beginning
at $i=i_0$ let $[i_0,i_1)$ be the $k$-interval and then take succesive
$k$-intervals $[i_1,i_2),\ldots,[i_{u-1},i_u)$ until reaching a $k$-interval
$[i_u,i_{u+1})$ with transient $k$-value. Then $[i,i_{u+1})$ is 
the $k+1$-interval of $i$. (Incidentally, suppose $U(i)$. Then $[i,i+1)$ 
is the 1-interval of $i$ which is transient.  But then $[i,i+1)$ is the 
2-interval of $i$ and is transient.  For all $k$ $[i,i+1)$ is the 
$k$-interval of $i$ and is transient.)
Let $x_1,\ldots,x_u,y_{u+1}$ be the succesive $k$-values of the intervals.
Let $\alpha$ be the equivalence class of $x_1\cdots x_u$ in $\Sigma P_k$.
Then $i$ has $k+1$-value $\beta=\alpha y_{u+1}$.  This value
is persistent if $\alpha$ is persistent and transient if $\alpha$ is
transient.  This defines $E_{k+1},P_{k+1},T_{k+1}$, completing the
induction. Our construction has assured that the $k$-value of $i$ 
determines the Ehrenfeucht value of the $k$-interval of $i$, though it
may have even more information.  
\par Now let us fix $i$ and look at the distribution of its $k$-value $V^k$.
We show, by induction on $k$, that for every persistent $\beta$
$ \Pr[V^k=\beta]=c_{\beta}+o(1)$ while for every transient $\beta$
$\Pr[V^k=\beta]=(c_{\beta}+o(1))p$.  Here each $c_{\beta}$ is a 
positive constant.  Assume the result for $k$ and
set $p_{\beta}=\Pr[V^k=\beta]$ for all $\beta\in E_k$. Let $p^*$ be the
probability that $V^k$ is transient so that $p^*\sim cp$, $c$ a 
positive constant.  Let $x_1,\ldots,x_u,y_{u+1}$ be the successive
$k$-values of the $k$-intervals beginning at $i$, stopping at
the first transient value.  We can assume these values are taken
independently from the inductively defined distribution on $E_k$.
The distribution of the first transient value is the conditional
distribution of $V^k$ given that $V^k$ is transient so the 
probability that it is some transient $y$ is $d_y+o(1)$
where $d_y=c_y/\sum c_{y'}$, the sum over all transient 
$y'$. Note all $d_y$ are positive constants.   
\par The key to the argument is the distribution for the
Ehrenfeucht equivalence
class $\alpha$ for the finite sequence $x_1\cdots x_u\in \Sigma P_k$.
Let $M$ be the set of all such equivalence classes.
Let $L_u$ be the event that precisely $u$ persistent $x$'s are
found and then a transient $y$.  Then $\Pr[L_u]=(1-p^*)^up^*$
precisely. For $\beta\in P_k$ let $p_{\beta}^+=p_{\beta}/(1-p^*)$,
the conditional probability that $V^k=\beta$ given that $V^k$ is
persistent.  Note that  (as $p\rightarrow 0$)
\[ p_{\beta}^+\sim p_{\beta}\sim c_{\beta}  \]
Conditioning on $L_u$ the $x_1,\ldots,x_u$ are 
mutually independent with distributions given by the $p_{\beta}^+$.
Define on $M$ a {\em Markov Chain} (see \S 3.3)  
with transition probability $p_{\beta}^+$ from  and $\alpha$ to
$\alpha+\beta$. We let $M(p)$ denote this Markov Chain.  Observe
that the {\em set} of states $M$ is independent of $p$ and the
nonzeroness of the transition probabilities is independent of
$p\in (0,1)$ though the actual transition probabilities do depend
on $p$.  
There is a particular state $O$ representing the
null sequence.  Let $f(u,\alpha)$ be the probability of being at
state $\alpha$ at time $u$, beginning at $O$ at time zero.  Then
$f(u,\alpha)$ is precisely the conditional distribution for $\alpha$
given $L_u$.  But therefore, letting $W$ denote the Ehrenfeucht
equivalence class,
\[ \Pr[W=\alpha]=\sum_{u=0}^{\infty} f(u,\alpha)(1-p^*)^up^*  \]
Let $M^o$ be the Markov Chain on the same set with transition
probability $c_{\beta}$ from $\alpha$ to $\alpha+\beta$ and
let $f^o(u,\alpha)$ be the probability of going from $O$ to 
$\alpha$ in $u$ steps under $M^o$.  Observe that $M^o$ is the
limit of $M(p)$ as $p\rightarrow 0$ in that taking the
limit of any (1-step) transition probability in $M(p)$ as
$p\rightarrow\infty$ gives the transition probability in 
$M^o$.
\par Now we need some Markov Chain asymptotics. Assume $\alpha$ is 
transient.
We claim (recall $p^*\sim cp$)
\[ \Pr[W=\alpha]\sim\left[c\sum_{u=0}^{\infty}f^o(u,\alpha)\right]p \]
and that the interior sum converges. In general the probability of
remaining in a transient state drops exponentially in $u$ so there 
exist constants $K,\epsilon$ so that $f^o(u,\alpha)<K(1-\epsilon)^u$
for all $u$ giving the convergence.  Moreover there exists $\epsilon_1,
\epsilon_2,K_1$ so that for all $0<p<\epsilon_1$ we bound uniformly
$f(u,\alpha)<K_1(1-\epsilon_2)^u$ for all $u$. Pick $\epsilon_3\leq
\epsilon_1$ so that for $0<p<\epsilon_3$ we have $p^*\leq 2cp$.
 For any positive $\delta$
we find $U$ so that for $0<p<\epsilon_3$
\[ \frac{\sum_{u=U}^{\infty} f(u,\alpha)(1-p^*)^up^*}{p}<
   \sum_{u=U}^{\infty} K_1(1-\epsilon_2)^u\frac{p^*}{p}
   \leq \frac{2cK_1}{\epsilon_2}(1-\epsilon_2)^U< \frac{\delta}{2} \]
For any fixed $0\leq u<U$ we have $\lim_{p\rightarrow 0}f(u,\alpha)=
f^o(u,\alpha)$ so that
\[ \lim_{p\rightarrow 0}\frac{\sum_{0\leq u<U} f(u,\alpha)
(1-p^*)^up^*}{p} = \sum_{0\leq u<U}cf^o(u,\alpha)  \]
With $U$ sufficiently large this may be made within $\delta/2$ of
$c\sum_0^{\infty}f^o(u,\alpha)$.  But this holds for $\delta$
arbitrarily small, giving the claimed asymptotics of $\Pr[W=\alpha]$.

\vspace{1cm}

{\em Remark.} The rough notion here is that the probability of
having a transient $k+1$-value is dominated by having few persistent
$k$-intervals and then a transient $k$-interval.  The transient
$2$-intervals all had at most $s$ persistent $1$-intervals. The
situation changes with $3$-intervals.  Recall $Ba_i$ consisted
of at least $s$ ones each preceeded by at least $s$ zeroes and
then two ones $i$ apart.  Consider an arbitrarily long
grouping of 2-intervals of 2-value $Ba_i$ but, say, with none
of the form $Ba_3$, i.e., $1001$ not appearing and then, say,
follow the last one, say $Ba_1$, with a one so that the $3$-
interval ends $111$.  For every $u$ there is a $\sim c_up$
probability of this being the $3$-interval with $u$ such
$2$-intervals and $c_u>0$ but all such $3$-intervals would
be considered transient since a persistent sequence in $\Sigma P_2$
must surely contain every value in $P_2$.

\vspace{1cm}

\par Now suppose $\alpha$ is persistent.  Again we have the
precise formula
\[    \Pr[W=\alpha]=\sum_{u=0}^{\infty}f(u,\alpha)(1-p^*)^up^* \]
only this time it is the tail of the sum that dominates.  As
$\alpha$ is persistent there is a limiting probability
$L=\lim_{u\rightarrow\infty}f^o(u,\alpha)$ with $L>0$ and
furthermore the $M(p)$ approach $M^o$ in the sense that
\[ L=\lim_{p\rightarrow 0}\lim_{u\rightarrow\infty}f(u,\alpha)  \]
We claim
\[ \Pr[W=\alpha]=L+o(1)  \]
For any $\delta>0$ there exist $\epsilon$ and $U$ so that for $p\leq\epsilon$
and $u\geq U$ we have
\[ L-\delta < f(u,\alpha) < L+\delta  \]
Then, as $\sum_{u=0}^{\infty} L(1-p^*)^up^* = L$,
\[ |\Pr[W=\alpha]-L| \leq \delta\sum_{u=U}^{\infty}(1-p^*)^up^* +
(L+1)\sum_{0\leq u<U}(1-p^*)^up^*  \]
For fixed $U$ the second sum is $o(1)$ (as $p^*\rightarrow 0$) while the
first sum is less than $\delta$ so the entire expression is less than 
$2\delta$ for $p$ sufficiently small. As $\delta$ was arbitrary
this gives the claim.
\par Recall that the $k+1$-value of the full $k+1$-interval is 
a pair consisting of the Ehrenfeucht value $W$ just discussed
and the $k$-value of the first transient type $y_{u+1}$.
The transient type's value has a limiting distribution which is
independent of $W$, for conditional on any $L_u$ the distribution
on $y_{u+1}$ is the same.  All possible $y\in T_k$ have a limiting
probability $d_y\in (0,1)$.  Hence the probability of a $k+1$-value
being $\beta=\alpha y$ is simply the product of the probabilities and
hence approaches a constant if $\alpha$, and hence $\beta$, is
persistent and is $\sim cp$ if $\alpha$, and hence $\beta$ is
transient.  This completes the inductive argument for the
limiting probabilities of the $k$-values of the $k$-intervals.

\par We now let $L=L^k$ be the length of the $k$-interval of $i$
and find bounds on the distribution of $L$.  A simple induction
shows that if the sequence $1\cdots 1$ of $k$ ones appears after
$i$ then the $k$-interval of $i$ ends with this sequence or
possibly before.  Thus we get the crude bound
\[ \Pr[L>ka]<(1-p^k)^a  \]
so that asymptotically
\[ Pr[L>\alpha p^{-k}]<e^{-c\alpha}  \]
where $c$ is a positive constant.  In fact, this gives the correct order
of magnitude, $L$ is (speaking roughly) almost always on the order of 
$p^{-k}$. We claim that there are positive constants $\epsilon_t,c_t$
so that 
\[ \Pr[L^t>\epsilon_tp^{-t}]>c_t \] 
The argument is by induction,
for $t=1$ the random variable $L^1$ is simply the number of trials until
a success which occurs with probability $p$ and the distribution is
easily computable.  Assume this true for $t$ and let (as previously
shown) $e_tp$ be the asymptotic probability that a $t$-interval will
be transient.  Pick $f_t$ positive with $f_te_t<.5$.
With probability at least $.5$, the first
$\gamma=f_tp^{-1}$ $t$-intervals after $i$ will be persistent.  Conditioning
on an interval being persistent is conditioning on an event that 
holds with probability $1-o(1)$ so that each of these $t$-intervals
will have length at least $\epsilon_tp^{-t}$ with probability at
least $c_t-o(1)$.  As the lengths are independent with conditional
probability at least $.99$ at least $c_t\gamma/2$ of the intervals
have length at least $\epsilon_tp^{-t}$. Thus with probability at
least, say $.4$ the total length $L^{t+1}$ is at least $c_t\gamma
\epsilon_tp^{-t}/2$ which is $\epsilon_{t+1}p^{-(t+1)}$ for an
appropriate constant $\epsilon_{t+1}$, completing the induction.
\par Up to now the relation between $p$ and $n$, the number of 
integers, has not appeared.  Recall that $p\rightarrow 0$ and
$n\rightarrow\infty$ so that $np^k\rightarrow\infty$ but
$np^{k+1}\rightarrow 0$.
Now begin at $i=i_0=1$ and generate
the $k$-interval $[i_0,i_1)$. Then generate the $k$-interval
$[i_1,i_2)$ beginning at $i_1$ and continue. 
(We do this with $k$ fixed.  Even if one of the intervals
is transient we simply continue with $k$-intervals.
Again we imagine continuing forever through the integers.)
Let $N$ be that
maximal $u$ for which $i_u-1\leq n$, so that we have split
$[n]$ into $N$ $k$-intervals plus some excess.  As each
sequence of $k$ ones definitely will end a $k$-interval $N$
is at least the number of disjoint subintervals of $k$ ones.
Simple expectation and variance calculations show that
$N>.99np^k$ almost surely.  On the other side set,
with  foresight, $c=4c_k^{-1}\epsilon_k^{-1}$ . 
If $N<cnp^k$ then the sum of the lengths
of the first $cnp^k$ $k$-intervals would be less than $n$.
But these lengths are independent identically distributed 
variables and each length is at least $\epsilon_kp^{-k}$
with probability at least $c_k$ so that almost surely at
least $c_kcnp^k/2$ of them would have length at least
$\epsilon_kp^{-k}$ and thus their total length would be
at least $(cc_k\epsilon_k/2)n>n$.  That is, almost surely
\[ C_1np^k < N < C_2np^k  \]
where $C_1,C_2$ are absolute constants.  
\par Let $\beta_1,\ldots,\beta_N$ be the $k$-values of the 
$k$-intervals generated by this procedure.
Now we make two claims about this procedure.  We first claim
that almost surely none of the $\beta_i$ are transient.
Each $\beta_i$ has probability $\sim cp$ of being transient
so the probability that some $\beta_i$, $1\leq i\leq C_2np^k$
is transient is at most $\sim (cp)C_2np^k=\Theta(np^{k+1})=o(1)$.
And almost surely $N<  C_2np^k$, proving the claim.
\par Let $A_1\cdots A_R$
be any {\em fixed} sequence of
elements of $P_k$.  The second claim is that almost surely
$A_1\cdots A_r$ appears as a subsequence of the $\beta$ sequence,
more precisely that almost surely there exists $i$ with 
$1\leq i\leq N-R$ so that $\beta_{i+j}=A_j$ for $1\leq j\leq R$.
(For technical reasons we want the subsequence not to start with
$\beta_1$.)  As each $\beta_i$ has a positive probability of 
being any particular $x\in P_k$ and the $\beta_i$ are independent
and $C_1np^k\rightarrow\infty$ almost surely this fixed sequence
will appear in the first $C_1np^k$ $\beta$'s.  And almost surely
$N>C_1np^k$, proving the claim.  
\par We have a third claim that is somewhat technical.  For any
$1\leq j\leq k$ let $\beta_1,\ldots,\beta_u$ denote the $j$-values
of the successive $j$-intervals starting at one, where $\beta_u$
is the last such interval that is in $P_j$.  We know that almost
surely $\beta_1\cdots\beta_u$ is persistent in $\Sigma P_j$. 
We claim further that almost surely $\beta_2\beta_3\cdots\beta_u$ is
persistent in $\Sigma P_j$.  It suffices to show this for any
particular $j$ as there are only a finite number of them.  For
any integer $A$ we have $u-1\geq A$ almost surely and the
probability that $\beta_2\cdots\beta_{A+1}$ is transient goes
to zero with $A$ so almost surely $\beta_2\cdots\beta_u$ is
persistent. Let us call $[b,c)$ a super $k$-interval (for a 
given $U$) if it is a $k$-interval and further for every $1\leq
j\leq k$ letting $\beta_1,\ldots,\beta_u$ be the successive
$j$-values of the $j$-intervals beginning at $b$ and stopping
with the last persistent value - that then $\beta_2\beta_3\cdots
\beta_u$ is persistent in $\Sigma P_j$.  So almost surely
the $k$-interval $[1,i_1)$ is a super $k$-interval.
\par We shall show, for an appropriate sequence $A_1,\ldots,A_R$,
that all $U$ satisfying the above three claims give models
$<n,C,U>$ which have the same Ehrenfeucht value.
\par We first need some glue.  Call $[a,b)$ an incomplete
$k$-interval (with respect to some fixed arbitrary $U$) if
the $k$-interval beginning at $a$ is not completed by $b-1$.
Suppose $[a,b)$ is an incomplete 
$k$-interval and $[b,c)$ is a persistent super
$k$-interval. We claim $[a,c)$ is a persistent $k$-interval.  
The argument is by induction on $k$.  For $k=1$, $[a,b)$ must
consist of just zeroes while $[b,c)$ consists of at least $s$
zeroes followed by a one.  But then so does $[a,c)$.  Assume
the result for $k$ and let $[a,b)$ be an incomplete $k+1$-interval
and $[b,c)$ be a persistent $k+1$-interval.  We split $[a,b)$
into a (possibly empty) sequence $x_1,x_2,\ldots,x_r$ of
persistent $k$-intervals followed by (possibly null) incomplete
$k$-interval $[a^+,b)$ with value, say, $y$.
We split $[b,c)$ (renumbering for convenience)
into a sequence $x_{r+1},\ldots,x_s,y_{s+1}$ of $k$-intervals,
all persistent except the last which is transient.   Then,
by induction, $y+x_{r+1}$ is a persistent $k$-interval with
some value $x_{r+1}'$.  Then $[a,c)$ splits into $k$-intervals
with values $x_1,\ldots,x_r,x_{r+1}',x_{r+2},\ldots,x_s,y_{s+1}$.
By the {\em super-}persistency $x_{r+2}\cdots x_s$ is persistent
in $\Sigma P_k$ and hence (see \S 3.2) so is $x_1\cdots x_rx_{r+1}'
x_{r+2} \cdots x_s$ and therefore $[a,c)$ is a {\em persistent} 
$k+1$- interval.
\par Now let $<[n],C,U>$ be any model that meets the three 
claims above, all of which hold almost surely for $p$ in
this range.  We set $i=i_0=1$ and find successive $k$-intervals
$[i_0,i_1),[i_1,i_2),\ldots$ until $[i_{u-1},i_u)$ and then
$U$ on $[i_u,n]$ gives an incomplete $k$-interval.  By the
third claim $[1,i_1)$ is superpersistent and so the
``interval'' $[i_u,n]\cup[1,i_1)$ (going around the corner)
is $k$-persistent.  Hence we have split $[n]$ (now thinking
of it as a cycle with $1$ following $n$ ) into $k$-persistent
intervals with $k$-values $x_1,x_2,\ldots,x_u$.  The $k$-value
for $x_1$ may be different from that for $[1,i_1)$ but the
others have remained the same. This sequence contains
the sequence $A_1\cdots A_R$ described in \S 3.5 .
But this implies (see \S 3.6) that the Ehrenfeuct value is determined,
completing the proof.

\section{Background.}
\subsection{The Ehrenfeucht Game.}
Let $A$ be a fixed finite alphabet (in application $A$ is $P_k$ or $\{
0,1\}$) and
$t$ a fixed positive integer.  We consider the space $\Sigma A$ of
finite sequences $a_1\cdots a_u$ of elements of $A$. We can associate
with each sequence a model $<[u],\leq,f>$ where $f:[u]\rightarrow A$
is given by $f(i)=a_i$.  For completeness
we describe the $t$-round Ehrenfeucht Games on sequences $a_1\cdots a_u$
and $a_1'\cdots a_{u'}'$.  There are two players, Spoiler and Duplicator.
On each round the Spoiler first selects one term from either sequencs and
then the Duplicator chooses a term from the other sequence.  Let $i_1,
\ldots i_t$ be the indices of the terms chosen from the first sequence,
$i_q$ in the $q$-th round and let $i_1',\ldots i_t'$ denote the corresponding
indices in the second sequence. For Duplicator to  win he must first assure
that $a_{i_q}=a_{i_q'}'$ for each $q$, i.e. that he selects each round
the same letter as Spoiler did.  Second he must assure that for all 
$a,b$ 
\[ i_a<i_b\Leftrightarrow i_a'<i_b' \mbox{ and } i_a=i_b\Leftrightarrow
i_a'=i_b' \]
(It is a foolish strategy for Spoiler to pick an already selected term
since Duplicator will simply pick its already selected counterpart but
this possiblity comes in in the Recursion discussed later.)  This is a 
perfect information game so some player will win.  Two sequences are called
equivalent if Duplicator wins.  Ehrenfeucht showed that this is an equivalence
class and that two sequences are equivalent if their models 
have the same truth value on all sentences of quantifier depth at most $t$.
We let $M$ denote the set of equivalence classes which is known to be a finite
set.  $\Sigma A$ forms a semigroup
under concatenation, denoted $+$, and this operation filters
to an operation, also denoted $+$, on $M$.  We use $x,y,\ldots$
to denote elements of $M$: $x+y$ their sum; $O$ is the equivalence
class of the null sequence which acts as identity.  We associate
$a\in A$ with the sequence $a$ of length one and its equivalence
class (which contains only it),, also called $a$.  We let $jx$
denote $x+\ldots+x$ with $j$ summands.  From analysis of the 
Ehrenfeucht game (see \S 3.4)
it is known that there exists $s$ (for definiteness
we may take $s=3^t$) so that:
\[  jx=kx \mbox{ for all } j,k\geq s,\mbox{  } x\in M \]

\par {\em Example.} With $A=\{0,1\}$ we naturally associate
sequences such as $101$ with $<\{1,2,3\},\leq,f>$ with
$f(1)=1,f(2)=0,f(3)=1$.  The addition of $101$ and $1101$
is their concatenation (in that order) $1011101$.  The first
order language has as atomic formulas $x\leq y$, $x=y$ and
$f(x)=a$ for each $a\in A$.  The sentence
\[ \exists_x\exists_y\exists_z f(x)=1\wedge f(y)=0\wedge f(z)=1
\wedge x<y \wedge y<z  \]
is satisfied by $01110001$ but not by $000111000$ so these
are in different equivalence classes with $t=3$. We could also
write that $101$ appears as consecutive terms with
\[ \exists_x\exists_y\exists_z f(x)=1\wedge f(y)=0\wedge f(z)=1\wedge
x<y\wedge y<z\wedge \neg\exists_w[(x<w\wedge w<y)\vee (y<w\wedge w<z)]
\]
Informally we would just say $\exists_x f(x)=f(x+1)=f(x+2)=1$ but
the quantifier depth is four. 

\subsection{Persistent and Transient.}
\par Definition and Theorem.  We call $x\in M$ {\em persistent} if
\[ \forall_y\exists_z x+y+z=x \hspace*{2cm}(1)  \]
\[ \forall_y\exists_z z+y+x=x \hspace*{2cm}(2)  \]
\[ \exists_p\exists_s\forall_y p+y+s=x \hspace*{2cm} (3) \]
These three properties are equivalent.  We call $x$ {\em transient} if
it is not persistent.
\\ Proof of Equivalence.
\par $(3)\Rightarrow(1):$ Take $z=s$, regardless of $y$.  Then
\[ x+y+z=(p+y+s)+y+s=p+(y+s+y)+s=x  \]
\par $(1)\Rightarrow(3):$ Let $R_x=\{x+v:v\in M\}$.
We first claim there exists $u\in M$ with $|R_x+u|=1$, i.e., all $x+y+u$
the same.  Otherwise take $u\in M$ with $|R_x+u|$ minimal and say
$v,w\in R_x+u$.  As $R_x+u\subseteq R_x$ we write $v=x+u_1, w=x+u_2$.
\relax From$(1)$, with $y=u_1$, we have $x=v+u_3$ and thus $w=v+u_4$ with
$u_4=u_3+u_2$.  Then
\[ w+su_4=v+(s+1)u_4=v+su_4  \]
Adding $su_4$ to $R+u$ sends $v,w$ to the same element so $|R+u+su_4|<
|R+u|$, contradicting the minimality.  Now say $R_x+u=\{u_5\}$.   Again
by $(1)$ there exists $u_6$ with $u_5+u_6=x$.  Then $R_x+(u+u_6)=\{x\}$
so that $(3)$ holds with $p=x, s=u+u_6$. 
\par By reversing addition (noting that $(3)$ is selfdual while the
dual of $(1)$ is $(2)$) these arguments give that $(3)$ and $(2)$ are
equivalent, completing the proof.
\par Let $x$ be persistent and consider $v=x+y$.  Let $z$ be such that
$x+w+z=x$ for all $w$.  Then for all $w$ $v+w+(z+y)=(x+(y+w)+z)+y=x+y=v$
and hence $v$ is persistent.  Dually,if $x$ is persistent $y+x$ is
persistent.  Together
\[ \mbox{If }x\mbox{ is persistent then } w_1+x+w_2\mbox{ is persistent} \]
for any $w_1,w_2\in M$.
\par From $(1)$ the relation $x\equiv_Ru$ defined by $\exists_v(x+v=u)$
is an equivalence relation on the set of persistent $x\in M$. We let
$R_x$ denote the $\equiv_R$-class containing $x$ so that
\[ R_x=\{x+v:v\in M\}  \]
\relax From$(2)$ the relation $x\equiv_Lu$ defined by $\exists_v(v+x=u)$
is also an equivalence relation on the set of persistent $x\in M$.
We let $L_x$ denote the $\equiv_L$-class containing $x$ so that
$L_x=\{v+x:v\in M\}$.  Let $x$ be persistent and let $p,s$ (by $(3)$)
be such that $p+z+s=x$ for all $x$.  Setting $z=O$, $x=p+s$. Thus
for all $z$
\[ x+z+x=(p+s)+z+(p+s)=p+(s+z+p)+s=x  \]
\par Let $R_x,L_y$ be equivalence classes under $\equiv_R,\equiv_L$
respectively.   Then $x+y\in R_x\cap L_y$.  Let $z\in R_x\cap L_y$.
Then there exist $a,b$ with $x=z+a$ and $y=b+z$ so that 
$x+y=z+(a+b)+z$.  But as $z$ is persistent the  above argument (with
$z$ as $x$ and $a+b$ as $z$) gives $z+(a+b)+z=z$.  Thus
\[ R_x\cap L_y=\{x+y\}\mbox{ for all persistent }x,y  \]

\vspace{1cm}

\par {\em Remarks.}  Let $A=\{0.1\}$.  A sequence $a_1\cdots a_u$ is
transient if and only if there is a sentence $Q$ of quantifier depth
at most $t$ so that $a_1\cdots a_u$ fails $Q$ but there is an extension
to $a_1\cdots a_ua_{u+1}\cdots a_v$ which satisfies $Q$ such that all
further extensions $a_1\cdots a_va_{v+1}\cdots a_w$ also satisfy $Q$.
For example, with $t=4$, let $Q$ be the existence of a block $101$.
If a sequence does not satisfy $Q$ then the extension given by
adding $101$ does satisfy $Q$ and all further extensions will satisfy
$Q$.  Thus for $a_1\cdots a_u$ to be persistent for $t=4$ it must
contain $101$ and indeed all blocks of length three.  We think of
property $(3)$ of persistency as indicating that a persistent
sequence is characterized by $p$, its prefix, and $s$, its suffix.
There are properties such as $\exists_x f(x)=1\wedge\neg\exists_yy<x$
that depend on the left side of the sequence, in this case the value
$f(1)$.  There are other properties such as $\exists_x f(x)=1\wedge
\neg\exists_ yx<y$ which depend on the right side of the sequence.
There will be sequences with values $p,s$ for the left and right
side respectively so that the Ehrenfeucht value of the sequence
is now determined, regardless of what is placed in the middle.

\vspace{1cm}

\par {\em Remarks.} Certain sentences  $Q$ have the property that
if any $a_1\cdots a_u$ satisfies $Q$ then all extensions $a_1\cdots
a_ua_{u+1}\cdots a_v$ satisfy $Q$.  The sentence that the first
term of the sequence is 1 has this property; the sentence that the
last term of the sequence is 1 does not have this property.  Call
such properties unrighteous, as they (roughly) do not depend on
the right hand side of the sequence.  Sequences with Ehrenfeucht 
value in a given $R_x$ all have the same truth value for all
unrighteous properties.  Sequences with Ehrenfeucht value in a 
given $L_x$ would all have the same truth value for all 
(correspondingly defined) unleftuous properties.
 
\subsection{The Markov Chain.}
Now consider a probability distribution over
$A$, selecting each $a$ with nonzero probability $p_a$.  This naturally
induces a distribution over $A^u$, the sequences of length $u$,
assuming each element is chosen independently.  This then leads to a 
distribution over the equivalence classes $M$.  For all $u\geq 0$, $x\in M$
let $P_u(x)$ be the probability that a random string $a_1\cdots a_u$ is
in class $x$.  On $M$ we define a Markov Chain, for each $x$ the transition
probability from $x$ to $x+a$ being $p_a$.  
\par In Markov Chain theory the states $x\in M$ are split into persistent
and transient, a state $x$ is persistent if and only if it lies in a 
minimal closed set.  We claim Markov Chain persistency is precisely
persistency as defined by $(1),(2),(3)$.  If $C$ is closed and $x\in C$
then $R_x\subseteq C$ and $R_x$ is itself closed. If $x$ satisfies
$(1)$ then $R_u=R_x$ for all $u=x+y\in R_x$ so $x$ is Markov Chain
persistent.  Conversely if $x$ is Markov Chain persistent then $R_x$
must be minimal closed so $R_u=R_x$ for all $u=x+y\in R_x$ and so
$x$ satisfies $(1)$.
\par The Markov Chain $M$ restricted to a minimal closed set $R_x$
is aperiodic since $x+sa\in R_x$ and $(x+sa)+a=x+sa$.  Hence from
Markov Chain theory when $x$ is persistent
$ \lim_{u\rightarrow\infty} P_u(x) $ exists.
\par A random walk on $M$, beginning at $O$, will with probability
one eventually reach a minimal closed set $R_x$ and then it must
stay in $R_x$ forever.  Let $P[R_x]$ denote the probability that
$R_x$ is the closed state reached.

\subsection{Recursion.}
Again let $A$ be a finite alphabet, $M$ the
set of equivalence classes of $\Sigma A$ and now specify some
$B\subseteq M$.  As $B$ is also a finite set we can define 
equivalence classes (with respect to the same constant $t$) on
$\Sigma B$, let $M^+$ denote the set of such classes. Now
let $b_1\cdots b_u$ and $b_1\cdots b_{u'}'$ be equivalent 
sequences of $\Sigma B$.  We claim that 
\[ b_1+\ldots+b_u=b_1'+\ldots+b_{u'}'  \]
as elements of $M$.  Let $s_1,\ldots,s_u,s_1',\ldots,s_{u'}'$
be specific elements of $\Sigma A$ in the repective $b_i$ or
$b_i'$ classes.  It suffices to give a strategy for Duplicator
with models $s_1+\ldots+s_u$ and $s_1'+\ldots+s_{u'}'$.  Suppose
Spoiler picks an element $x$ in, say, some $s_i$.  In the game
on $\Sigma B$ we know Duplicator has a winning reply to $b_i$
of some $b_{i'}'$.  Now Duplicator will pick some $x'$ in
$s_{i'}'$.  To decide the appropriate $x'$ in $s_{i'}'$ to
pick Duplicator considers a subgame on $s_i$ and $s_{i'}'$.  As
these are equivalent Duplicator will be able to find such $x'$
for the at most $t$ times that he is required to.
\par This general recursion includes the previous statement
that for all $j,k\geq s$ and any $x\in M$ we have $jx=kx$.
Here $B=\{x\}$ and this says that Duplicator can win the
$t$-move Ehrenfeucht game between a sequence of $j$ $x$'s and
a sequence of $k$ $x$'s; that is, that $<[j],\leq>$ and
$<[k],\leq>$ are equivalent - a basic result on Ehrenfeucht games.
In our argument we apply it inductively with $A=P_k$.  We know,
inductively, that all $k$-intervals having the same $k$-value
$x\in P_k$ have the same Ehrenfeucht value.  Now the $k+1$-interval
of $i$ is associated with
a sequence $x_1\cdots x_u\in \Sigma P_k$  and a ``tail'' $y_{u+1}\in T_k$.
We call two such $k+1$-intervals equivalent if $x_1\cdots x_{u'}'$ and
$x_1'\cdots x_{u'}'$ are equivalent in $\Sigma P_k$ and $y_{u+1}=
y_{u'+1}'$.  Now $x_1+\ldots+x_u=x_1'+\ldots+x_{u'}'$ and so
the $k+1$-intervals have equal Ehrenfeucht value.

\subsection{Cycles.}
Again let $M$ be the set of equivalence classes
on $\Sigma A$. Now consider cycles $a_1\cdots a_u$ (thinking of
$a_1$ following $a_u$) with $a_i\in A$ and consider equivalence classes
under the $(t+1)$-move Ehrenfeucht game.  Here we must preserve the
ternary clockwise predicate $C(x,y,z)$.  Any first move $a_i$ reduces
the cycle to a linear structure $a_i\cdots a_ua_1\cdots a_{i-1}$ of
the form $<[u],\leq,f>$ with an Ehrenfeucht value $x=x_i$.  Two
cycles are equivalent if they yield the same set of values $x_i\in M$.
\par For every persistent $x\in M$ let (by $(3)$) $p=p_x, q=s_x$ be
such that $x=p_x+y+s_x$ for all $y\in M$.  Let $P_x$ and $S_x$ be
fixed sequences (i.e., elements of $\Sigma A$) for these equivalence
classes and let $R_x$ be the sequence consisting of $S_x$ {\em in reverse
order} followed by $P_x$.  If the cycle $a_1\cdots a_u$ contains
$R_x$ as a subsequence then selecting $a_i$ as the first element
of $P_x$ gives a linear structure beginning with $P_x$ and ending
with $S_x$, hence of value $p_x+y+s_x=x$.  
\par Let $R\in \Sigma A$ be a specific sequence given by the 
concatenation of the above $R_x$ for all persistent $x\in M$.
Then we claim $R$ is a {\em universal sequence} in the sense 
that all $a_1\cdots a_u\in \Sigma A$ (for any $u$)
that contain $R$ as a subsequence  are equivalent.  For every
persistent $x\in M$ there is an $a_i$ so that $a_i\cdots a_{i-1}$
has value $x$.  Conversely every $a_i$ belongs to at most one of
the $R_x$ creating $R$ (maybe none if $a_i$ isn't part of $R$)
and so there will be an $R_x$ not containing that $a_i$.  Then
in $a_i\cdots a_{i-1}$ the subsequence $R_x$ will appear as an
interval.  Hence the value of $a_i\cdots a_{i-1}$ can be written
$w_1+x+w_2$, which is persistent.  That is, the values of
$a_i\cdots a_{i-1}$ are precisely the persistent $x$ and hence
the class of $a_1\cdots a_u$ in the circular $t+1$-Ehrenfeucht game
is determined.  

\subsection{Recursion on Cycles.}
Again let $A$ be a finite alphabet, $M$ the set of equivalence classes
in $\Sigma A$ and specify some $B\subseteq M$.  Suppose a cycle 
$a_1\cdots a_u$ on $A$ may be decomposed into intervals $s_1,\ldots,
s_r$ with Ehrenfeucht values $b_1\cdots b_r$.  Then the Ehrenfeucht
value of the cycle $b_1\cdots b_r$ determines the Ehrenfeucht value
of $a_1\cdots a_u$. The argument is the same as for recursion on
intervals.  Let $a_1\cdots a_u$ and $a_1'\cdots a_{u'}'$ be decomposed
into $s_1\cdots s_r$ and $s_1'\cdots s_{r'}'$ with Ehrenfeucht values
$b_1\cdots b_r$ and $b_1\cdots b_{r'}'$.  Spoiler picks $x$ in some
$s_i$.  In the game on cycles over $B$ Duplicator can respond 
$b_{i'}'$ to $b_i$.  Then Duplicator picks an $x'\in s_{i'}'$ so that
he can win the subgame on $s_i$ and $s_{i'}'$.
\par We apply this is \S 2 with $A=\{0,1\}$ and $B=P_k$.  Here the
$\beta\in P_k$ may have more information than the Ehrenfeucht value
but this only helps Duplicator.

\section{The Linear Model.}  
We have already remarked in \S 1 that Zero-One Laws generally
do not hold for the linear model $<[n],\leq,U>$ and that P.
Dolan has characterized those $p=p(n)$ for which they do.
Our main object in this section is the following {\em convergence}
result.
\\ {\bf Theorem 2.}  Let $k$ be a positive integer, and $S$ a 
first order sentence.  Then there is a constant $c=c_{k,S}$ so
that for any $p=p(n)$ satisfying
\[ n^{-\frac{1}{k}}\ll p(n)\ll n^{-\frac{1}{k+1}}  \]
we have
\[ \lim_{n\rightarrow\infty}\Pr[U_{n,p}\models S]= c  \]
\par Again we shall fix the quantifier depth $t$ of $S$ and
consider Ehrenfeucht classes with respect to that $t$. For each
$\beta\in P_k$ let $c_{\beta}$ be the constant defined in \S 2
as the limiting probability that a $k$-interval has $k$-value 
$\beta$.  Let $M$ be the set of equivalence classes of $\Sigma P_k$,
a Markov Chain as defined in \S 3, and for each $\equiv_R$-class
$R_x$ let $P[R_x]$, as defined in \S 3, be the probability that 
a random sequence $\beta_1\beta_2\cdots$ eventually falls into
$R_x$.
\par In $<[n],\leq,U>$ let $\beta_1\cdots \beta_N$ denote the
sequence of $k$-values of the successive $k$-intervals, denoted 
$[1,i_1),[i_1,i_2),\ldots$, from 1.
\par Set, with foresight, $\delta=10^{-2}3^{-t}$.
\par We shall call $U$ on $[n]$ {\em right nice} if it satisfies
two conditions.  The first is simply that all the $\beta_1,\ldots,
\beta_N$ described above are persistent.  We know from \S 2 that
this holds almost surely.  The second will be a 
particular universality condition.  Let $A_1\cdots A_R$ be
a specific sequence in $\Sigma P_k$ with the property that
for every $R_x$ and $L_y$ there exists a $q$ so that
\[ A_1\cdots A_q\in L_y \mbox{ and } A_{q+1}\cdots A_R\in R_x  \]
(We can find such a sequence for a particular choice of $R_x$ and
$L_y$ by taking specific sequences in $\Sigma P_k$ in those 
classes and concatenating them.  The full sequence is achieved by
concatenting these sequences for all choices of $R_x$ and $L_y$.
Note that as some $A_1\cdots A_q\in L_y$ the full sequence is
persistent.)  The second condition is that inside any interval
$[x,x+\delta n]\subset [1,n]$ there exist $R$ consecutive
$k$-intervals $[i_L,i_{L+1}),\ldots,[i_{L+R},i_{L+R+1})$ whose
$k$-values are, in order, precisely $A_1,\ldots,A_R$.  We claim
this condition holds almost surely.  We can cover  $[1,n]$ with
a finite number of intervals $[y,y+\frac{\delta}{3}n]$ and it
suffices to show that almost always all of them contain such a 
sequence, so it suffices to show that a fixed $[y,y+\frac{\delta}
{3}n]$ has such a sequence.  Generating the $k$-intervals from 1
almost surely a $k$-interval ends after $y$ and before $y+
\frac{\delta}{6}n$.  Now we generate a random sequence $\beta_1\cdots$
on an interval of length $\frac{\delta}{6}n$.  But constants do not affect
the analysis of \S 2 and almost surely $A_1\cdots A_R$ appears.
\par Now on $<[n],\leq,U>$ define $U^r$ by $U^r(i)$ if and only if
$U(n+1-r)$. $U^r$ is the sequence $U$ in reverse order.  Call $U$
left nice if $U^r$ is right nice.  Call $U$ nice if it is right nice
and left nice.  As all four conditions hold almost surely, the random
$U_{n,p}$ is almost surely nice.
\par Let $U$ be nice and let $\beta_1\cdots\beta_N$ and $\beta_1^r
\cdots \beta_{N^r}^r$ denote the sequences of $k$-values for 
$U$ and $U^r$ respectively and let $R_x$ and $R_{x^r}$ denote
their $\equiv_R$-classes respectively.  (Both exist since the
sequences are persistent.)
\\ {\em Claim.}  The values $R_x$ and $R_{x^r}$ determine the 
Ehrenfeucht value of nice $U$.
\par We first show that Theorem 2 will follow from the Claim.  Let
$R_x,R_{x^r}$ be any two $\equiv_R$-classes.  Let $U$ be random and
consider $<[\delta n],\leq,U>$.  The sequence of $k$-values lies
in $R_x$ with probability $P[R_x]+o(1)$.  The same holds for
$U^r$  on $[\delta n]$.  But $U^r$ examines $U$ on $[(1-\delta)n,n]$
so as $\delta<.5$ the values of the $\equiv_R$-classes are 
{\em independent} and so the joint probability of the values being
$R_x$ and $R_{x^r}$ respectively is $P[R_x]P[R_{x^r}]+o(1)$.  Given
the Claim $<[n],\leq,U>$ would then have a value $v=v(R_x,R_{x^r})$.
As
\[ \sum P[R_x]P[R_{x^r}]=\sum P[R_x]\sum P[R_{x^r}] = 1\times 1=1 \]
this would give a limiting distribution for the Ehrenfeucht value
$v$ on $<[n],\leq,U>$.
\par Now for the claim. Fix two models $M=<[n],\leq,U>$ and 
$M'=<[n'],\leq,U'>$, both nice and both with the same values
$R_x,R_{x^r}$.  Consider the $t$-move Ehrenfeucht game.  For
the first move suppose Spoiler picks $m\in M$.  By symmetry
suppose $m\leq \frac{n}{2}$. Let $[i_{r-1},i_r)$ be one of 
the $k$-intervals with, say, $.51n\leq i_r\leq .52n$.  We allow
Duplicator a ``free'' move and have him select $i_r$.  Let
$\beta_1\cdots\beta_N$ and $\beta_1'\cdots\beta_{N'}'$ be the
sequences of $k$-values for $M$ and $M'$ respectively.  Let $z$
be the class of $\beta_1\cdots\beta_r$. Since $U$ is nice this
sequence already contains $A_1\cdots A_R$ and hence is persistent
so $z\in R_x$. Let $z'$ be the class of $\beta_{r+1}\cdots\beta_N$.
By the same argument $z'$ is persistent.  In $M'$ inside of, say,
$[.5n,.51n]$ we find the block $A_1\cdots A_R$.  By the universality
property we can split this block into a segment in $L_z$ and 
another in $R_{z'}$.  Adding more to the left or right doesn't
change the nature of this split. Thus there is an interval
$[i_{r'-1}',i_{r'}')$ so that $\beta_1'\cdots\beta_{r'}'\in L_z$
and $\beta_{r'+1}'\cdots\beta_{N'}'\in R_{z'}$.  Spoiler plays
$i_{r'}'$ in response to $i_r$.  
\par The class of $\beta_1\cdots\beta_r$ is $z$ and $z\in R_x$.
The class $z'$ of $\beta_1'\cdots\beta_{r'}'$ is in $L_z$ and
$R_x$.  As $z\in L_z\cap R_x$, $z=z'$.  Thus $[1,i_r)$ under $M$
and $[1,i_{r'}')$ under $M'$ have the same Ehrenfeucht value.
Thus Duplicator can respond successfully to the at most $t$
moves (including the initial move $m$) made in these intervals.
Thus Spoiler may as well play the remaining $t-1$ moves on
$M_1=<[i_r,n],\leq,U>$ and $M_1'=<[i_{r'}',n'],\leq,U'>$.
These intervals have lengths $n_1\geq \frac{n}{3}$ and $n_1'\geq
\frac{n'}{3}$ respectively.  But now $M$ and $M'$ are both nice 
with respect to $\delta_1=3\delta$ - the sequence $A_1\cdots A_R$
still appears inside every interval of length $\delta n\leq\delta_1n_1$
in $M$ and $\delta_1n_1'$ in $M'$.  Hence we can apply the same 
argument for the second move - for convenience still looking at
Ehrenfeucht values with respect to the $t$ move game.  After $t$
moves we still have nice $M_t,M_t'$ with respect to $\delta_t\leq 10^{-2}$
so the arguments are still valid.  But at the end of $t$ rounds Duplicator
has won.

\end{document}